\documentclass[12pt,a4paper]{article}
\usepackage{amssymb,amsmath}
\usepackage{latexsym}
\font\emailfont=cmtt10
\font\tenBbb=msbm10 scaled 1200
\font\sevenBbb=msbm7 scaled 1200
\font\fiveBbb=msbm5 scaled 1200
\newfam\Bbbfam
\textfont\Bbbfam=\tenBbb
\scriptfont\Bbbfam=\sevenBbb
\scriptscriptfont\Bbbfam=\fiveBbb
\def\Bbb{\fam\Bbbfam\tenBbb}
\global\tracingstats=0
\newtheorem{thrm}{Theorem}[section]
\newtheorem{lem}[thrm]{Lemma}
\newtheorem{prop}[thrm]{Proposition}
\newtheorem{cor}[thrm]{Corollary}

\newcommand{\aR}{{\Bbb R}}
\def\sph#1{{\Bbb S}^{#1}}
\def\littleoh#1{o(#1)}
\def\bigoh#1{{\rm O}(#1)}
\newcounter{dummy}
\newcommand{\numnow}{(\theequation)}
\newcommand{\numbefore}{\setcounter{dummy}{\theequation}\addtocounter{dummy}{-1}(\thedummy)}
\newcommand{\numbeforethat}{\setcounter{dummy}{\theequation}\addtocounter{dummy}{-2}(\thedummy)}
\def\eps{\varepsilon}
\def\phi{\varphi}
\def\norma#1{|\!|#1|\!|}
\def\8{\infty}

\def\ip#1#2{\langle #1,#2\rangle}
\def\bigskip{\vskip 12pt}
\def\medskip{\vskip 8pt}
\def\cube#1{B_{\8}^{#1}}
\def\face{F_{n-k+j}}

\def\thinsquare{\vcenter{\hrule height 0.3pt \hbox{\vrule width0.3pt height6pt
\kern6.5pt \vrule width 0.3pt}
\hrule height 0.3pt}}
\def\sofproof{\hfill$\thinsquare$}
\newcommand{\lechol}{\forall\,}
\newcommand{\gee}{{\Bbb G}}

\newcommand{\haflip}{\frac{2}{\pi}}
\newcommand{\ex}{\textbf{E}\,}
\newcommand{\proof}{\emph{Proof}\/.\quad} \author{Yossi Lonke\thanks{\ \
Supported in part by the US-Israel Binational Science Foundation. This work forms
part of a dissertation prepared under the direction of Professor Joram
Lindenstrauss, to whom I am greatly indebted for much valuable advice.}}
\title{On random sections of the cube}
\date{December 14, 1998}
\begin{document}
\maketitle
\begin{abstract}
Let $f(j,k,n)$ denote the expected number of $j$-faces of a random
$k$-section of the $n$-cube. A formula for $f(0,k,n)$ is presented, and for
$j\geq 1$, a
lower bound for $f(j,k,n)$ is derived, which implies a precise asymptotic formula
for $f(n-m,n-l,n)$ when $1\leq l<m$ are fixed integers and $n\to\8$.
\end{abstract}
\section{Introduction}

The principal object in this paper is the expected number of $j$-dimensional
faces (in short, $j$-faces) of a random $k$-dimensional central section (in
short, $k$-section) of the $n$-cube $\cube{n}=[-1,1]^n$ in $\aR^n$. We denote
this number by $f(j,k,n)$. The normalized rotation invariant measure on the set
$G_{n,k}$ of all $k$-dimensional subspaces of $\aR^n$ provides the probabilistic
framework. 

Section~$2$ contains a calculation of the expected number of vertices of a random
$k$- section
of the $n$-cube. The result is:
\begin{equation}\label{eq:myvertex}
f(0,k,n)=2^k{n\choose k}\sqrt{\frac{2k}{\pi}}\int_0^{\8}e^{-kt^2/2}
\gamma_{n-k}(t\cube{n-k})\,dt,
\end{equation}
where $\gamma_{n-k}$ denotes the $(n-k)$-dimensional Gaussian probability
measure. 

In $\S 3$ we derive a lower bound for $f(j,k,n)$ for every $1\leq j<k<n$. The
main result is:
$$
\frac{f(0,k-j,n)}{f(j,k,n)}<
\sqrt{\frac{2}{\pi}}\Big(\frac{j(k-j)}{n-k+j}\Big)^{1/2}
\int_0^{\8}e^{-\frac{j(k-j)}{n-k+j}t^2/2}\gamma_j(t\cube{j})\,dt.
$$
The lower bound for $f(j,k,n)$ derived from this inequality, combined with (\ref{eq:myvertex}), leads
in some cases to asymptotically best possible results.  For example, in $\S 3$ we deduce
from it the following asymptotic formula, for fixed
integers $1\leq l<m$:
\begin{equation}\label{eq:myasympt}
f(n-m,n-l,n)\backsim\frac{(2n)^{m-l}}{(m-l)!},\quad\textrm{as}\ \ n\to\8. 
\end{equation}
The
notation $a_n\sim b_n$ means: $a_n/b_n\to 1$ as $n\to\8$. (\ref{eq:myasympt}) can be interpreted as follows: the probability that a
random fixed-codimensional subspace of $\aR^n$ intersects a fixed-codimensional face of the
$n$-cube, tends to $1$ as $n\to\8$. The formula (\ref{eq:myasympt}) itself follows also
from the work of Affentranger and Schneider. (See remark~1 of section~3 below). In \cite{AffSch}, they found a formula for the expected
number $\ex(f_j(\Pi_kB_1^n))$ of $j$-faces of an orthogonal
projection of an $n$-polytope $P$ onto a $k$-dimensional random
subspace. Formula (5) of \cite{AffSch} reads as follows:
\begin{equation}\label{eq:AffSch}
\ex(f_j(\Pi_kP))=f_j(P)-2\sum_{s\geq 0}\sum_{F\in{\cal F}_j(P)}\sum_{G\in{\cal F}_{k+1+2s}(P)}\beta(F,G)\gamma(G,P).
\end{equation}
Here ${\cal F}_j(P)$ denotes the set of $k$-faces of $P$, and $f_j(P)=\hbox{card}\,{\cal F}_j(P)$.
 $\beta(F,G)$ denotes the internal angle (\cite{Grun}, p.~297) of the face $G$ at
its face $F$, and $\gamma(G,P)$ --- the external angle (\cite{Grun}, p.~308) of~$P$ at
its face~$G$. It is shown in \cite{AffSch} that (\ref{eq:AffSch}) implies that if $0\leq j<k$ are
fixed integers, then as $n\to\8$,
\begin{equation}\label{eq:BHasymp1}
\ex(f_j(\Pi_kT^n))\sim\frac{2^k}{\sqrt{k}}{k\choose j+1}\beta(T^j,T^{k-1})
(\pi\log n)^{(k-1)/2}.
\end{equation}
Here $T^n$ stands for the regular $n$-simplex. \par
In a very recent work,
\cite{BorHenk}, B\"or\"oczky, Jr. and Henk showed that (\ref{eq:AffSch}) implies
the same asymptotic formula (\ref{eq:BHasymp1}) also for
$\ex(f_j(\Pi_kB_1^n))$, where $B_1^n$ is the regular cross-polytope. In addition, they found an
asymptotic formula for the internal angles $\beta(T^j,T^{k-1})$, when
$k/j^2\to\8$. Therefore if $j$ is fixed, $k$ is much larger than $j^2$ and $n$
much larger than $k$, then explicit estimates for $\ex(f_j(\Pi_kB_1^n))$ are
available. See \cite{BorHenk} for more details. 
Explicit asymptotic
formulas for
$\ex(f_j(\Pi_kT^n))$, were established
independently
by Vershik and Sporyshev (\cite{VerSp}), when $j,k$ are both proportional to~$n$
and ${n\to\8}$.

A simple duality argument shows that
$$\ex(f_j(\Pi_kB_1^n))=f(k-j-1,k,n).$$
Choose $j=k-1$ in (\ref{eq:BHasymp1}). Applying the result for $B_1^n$, one has
\begin{equation}\label{eq:vertexasymp}
f(0,k,n)=\ex(f_{k-1}(\Pi_kB_1^n))\sim \frac{2^k}{\sqrt{k}}(\pi\log n)^{(k-1)/2},
\qquad\textrm{as $n\to\8$}.
\end{equation}
The last asymptotic formula follows also from (\ref{eq:myvertex}). In fact, if
$\{g_i\}_{i=1}^m$ are independent $N(0,1)$ (that is, with mean~$0$ and
variance~$1$) Gaussian variables then
$\gamma_m(t\cube{m})$ coincides with the probability of the event $\{\max_{1\leq
i\leq m}|g_i|\leq t\}$. This probabilistic interpretation allows a
straightforward evaluation of the asymptotic behavior of the integral in
(\ref{eq:myvertex}), when $k$ is fixed and $n\to\8$. 

Formula (\ref{eq:myvertex}) also yields information about $f(0,k,n)$ for $k$ not
necessarily fixed. For example, if $k=n-1$, then the integral in
(\ref{eq:myvertex}) can
be computed and the result is:
\begin{equation}\label{eq:myfacets}
f(0,n-1,n)=\frac{2^nn}{\pi}\arctan\frac{1}{\sqrt{n-1}}\sim\frac{2^n\sqrt{n}}{\pi}.
\end{equation}
Particular values of the last formula were computed numerically in
\cite{BorHenk}. (Table~$2$). For the expected number of vertices of random
sections of fixed co-dimension, we have the following inequality, which is a
consequence of~(\ref{eq:myvertex}). $$f(0,n-d,n)\geq {n\choose
d}2^n\left(\frac{1}{\pi}\arctan\frac{1}{\sqrt{n- d}}\right)^d,\qquad (d\geq 1).$$
Equality holds for $d=1$. 

To obtain a lower bound for $f(j,k,n)$, it turns out that it is useful to know an
estimate for the Gaussian measure of a cone generated by a section of a face of a
cube. In $\S 3$ we find such an estimate, by modifying K.~Ball's calculation of
the maximal volume of a cube--section, based on Brascamp-Lieb's inequality.
(\cite{Balcub}). 

Dvoretzky's theorem on almost Euclidean sections asserts that there exists a
function
$k(\eps,n)\geq 1$, tending to infinity as $n\to\8$ for each fixed $\eps>0$, such that if
$K$ is an $n$-dimensional centrally symmetric convex body (that is, a convex
compact set in
$\aR^n$ with non-empty interior, satisfying $K=-K$), and $\eps>0$, then for each
$1\leq k\leq k(\eps,n)$ there exists a $k$-dimensional subspace $X$, and a linear
automorphism $T$ of
$X$ for which
\begin{equation}\label{eq:BMdist}
X\cap B_2^n\subset T(X\cap K)\subset (1+\eps)(X\cap B_2^n), \end{equation}
where $B_2^n$ denotes the Euclidean unit ball. The proof of Dvoretzky's theorem
in \cite{FLM} shows that ${k(\eps,n)\geq c\eps^2|\log\eps|^{-1}\log n}$, for some
absolute constant $c>0$. That proof determined the best possible dependence
of~$k$ on~$n$. The dependence of~$k$
on~$\eps$ was improved by Gordon \cite{Gord}, who discovered another proof of
Dvoretzky's theorem
with $k(\eps,n)\geq c\eps^2\log n$. Both proofs are probabilistic; they show that
not only
there exist almost Euclidean sections, but actually most sections are such. More
precisely, if $X$ is a random subspace whose dimension does not exceed
$k(\eps,n)$, then the probability that the section $X\cap K$ is
$(1+\eps)$-Euclidean (common terminology for expressing that (\ref{eq:BMdist})
holds), tends to~$1$ as $n\to\8$. These facts motivate an investigation of the
random $f$-vector $\{f(j,k,n)\} _{j=0}^{k-1}$, especially since it is well known
that every $k$-dimensional symmetric polytope that has $2n$ facets is affinely
equivalent to a $k$-section of an $n$-cube. \section{Vertices}
Let $G_{n,k}$ denote the set of $k$-dimensional subspaces of $\aR^n$. We will
denote its normalized rotation invariant measure by ``Prob''. Recall that this
measure is related to the normalized Haar measure $H$ of the orthogonal group
${\rm O}(n)$ by the equality $$\textrm{Prob\,}\{X\in B\}=H\{g\in\textrm{O}(n):
g[e_i]_{i=1}^k\in B\},$$ where $B$ is a Borel subset of $G_{n,k}$ and
$[e_i]_{i=1}^k$ is the $k$-dimensional subspace spanned by the first $k$ unit
vectors in $\aR^n$. Fix $X\in G_{n,k}$. For each ${0\leq j\leq k-1}$, the set of
$j$-faces
of the polytope $X\cap\cube{n}$ coincides with the set of intersections of
$(n-k+j)$-faces of~$\cube{n}$ with~$X$. Every $(n-k+j)$-face of $\cube{n}$ has
the same probability to be intersected.
Therefore if one particular $(n-k+j)$-face $F_{n-k+j}$ is fixed, then the
expected number of $j$-faces of the section $X\cap\cube{n}$ is equal to:
\begin{displaymath}
2^{k-j}{n\choose k-j}\textrm{Prob\,}
\{X\cap F_{n-k+j}\neq\emptyset\}.
\end{displaymath}
Let $C(F_{n-k+j})$ denote the cone generated by $F_{n-k+j}$:
$$C(F_{n-k+j})=\bigcup_{x\in F_{n-k+j}}\{tx: t\geq 0\}.$$ Put
$C_1(F_{n-k+j})=C(F_{n-k+j})\cap\sph{n-1}$. For every subspace $X$, $$
X\cap\face\neq\emptyset\Longleftrightarrow(X\cap\sph{n-1})\cap
C_1(\face)\neq\emptyset. $$
For $n=0,1,\dots$ we denote by $\sigma_n$ the normalized rotation-invariant
measure on the unit-sphere $\sph{n}$ in $\aR^{n+1}$. The next lemma will prove
useful for dealing with intersections of subsets of the sphere with random
subspaces.
\begin{lem}
Let $l,m,n$ be positive integers satisfying $l+m\geq n-1$. Suppose that
$A\subset\sph{m}$ and $B\subset\sph{l}$ are Borel subsets. Then for
$p=l+m-n+1$,
\begin{equation}\label{eq:mainlem}
\int_{{\rm O}(n)}\sigma_p(gB\cap
A)\,dH(g)=\sigma_l(B)\,\sigma_m(A).
\end{equation}
\end{lem}
To prove the lemma one observes that for fixed $A$ (resp. $B$) the integral
defines an invariant
measure on $\sph{l}$ (resp. $\sph{m}$); the conclusion follows from that. 

Lemma~$2.1$ is now applied to $B=X\cap\sph{n-1}$, which we denote by $\sph{k-1}$,
and to~${A=C_1(\face)}$.
For $l=k-1$ and $m=n-k+j$ equality (\ref{eq:mainlem}) becomes:
\begin{equation}\label{eq:basic}
\int_{{\rm O}(n)}\sigma_j(g\sph{k-1}\cap A)\,dH(g)=\sigma_{n-k+j}(A).
\end{equation}

We are ready to compute the expected number of vertices. The Gaussian measure in
$\aR^m$ whose density is $(2\pi)^{-m/2}\exp(-\sum_1^mx_i^2/2)$ is denoted by
$\gamma_m$. \begin{prop}
The expected number of vertices of a random $k$-dimensional central section of
the $n$-cube is given by the formula
$$f(0,k,n)=2^k{n\choose k}\sqrt{\frac{2k}{\pi}}\int_0^{\8}e^{-kt^2/2}
\gamma_{n-k}(t\cube{n-k})\,dt.$$
\end{prop}
\proof
For each $g\in
\textrm{O}(n)$ we have
$$g\sph{k-1}\cap C_1(F_{n-k})=
(\textrm{span}(g\sph{k-1})\cap C(F_{n-k}))\cap\sph{n-1}.$$ For almost every $g$
the intersection
$\textrm{span}(g\sph{k-1})\cap C(F_{n-k})$ is either the origin itself, or else a
ray emanating from the origin. Therefore the intersection $g\sph{k-1}\cap
C_1(F_{n-k})$ is either empty or a singleton, for almost every~$g$. Choose $j=0$
in (\ref{eq:basic}), with $A=C_1(F_{n-k})$. Since the measure $\sigma_0$ is
concentrated on two points giving mass $1/2$ to each,
we deduce from $\numnow$ that
\begin{equation}\label{eq:Prob}
\textrm{Prob\,}\{X\cap
F_{n-k}\neq\emptyset\}=2\sigma_{n-k}(C_1(F_{n-k})). \end{equation}
To compute
the r.h.s of $\numnow$, consider an $(n-k)$-dimensional cube of edge-length $1$
inside
$\aR^{n-k+1}$, at a distance $\sqrt{k}$ from the origin, form the cone it
generates, and compute
the measure of its intersection with the sphere $\sph{n-k}$. Invoking polar
coordinates we see that
$$\sigma_{n-k}(C_1(F_{n-k}))=\gamma_{n-k+1}(C(F_{n-k})).$$ By rotational symmetry
of the Gaussian measure we may assume that $F_{n-k}$ is specifically
the set $\{x: |x_i|\leq1,\, 1\leq i\leq n-k,\,x_{n-k+1}=\sqrt{k}\}$. The
intersection of the hyper-plane $\{x_{n-k+1}=t\}$ with $C(F_{n-k})$ is an 
$(n-k)$-dimensional cube of edge-length $\frac{t}{\sqrt{k}}$. Therefore by
Fubini's theorem \begin{displaymath}
\aligned
\gamma_{n-k+1}(C(F_{n-k}))&=\frac{1}{\sqrt{2\pi}}\int_0^{\8}e^{-
t^2/2}\gamma_{n-k}(\frac{t}{\sqrt{k}}\cube{n-k})\,dt\cr
&=\sqrt{\frac{k}{2\pi}}\int_0^{\8}e^{-kt^2/2}\gamma_{n-k}(t\cube{n-k})\,dt.
\endaligned
\end{displaymath}
The last equality, together with $\numnow$, implies the desired formula.\sofproof
\bigskip
The next lemma points out the precise asymptotic behavior of $f(0,k,n)$ when $k$
is fixed and $n\to\8$, and also that of $f(n-m,n-l,n)$, when $l,m$ are fixed and
$n\to\8$. (To be used in~$\S 3$.)
\begin{lem}
Suppose that $\{\alpha_n\}_{n=1}^{\8}$ is a sequence of real numbers that has a
positive limit $\alpha$. Then as
$n\to\8$,
\begin{equation}\label{eq:myasymp}
\int_0^{\8}e^{-\alpha_n t^2/2}\gamma_n(t\cube{n})\,dt\sim
\Gamma(\alpha)\frac{\pi^{\alpha/2}}{\sqrt{2}}\frac{(\log n)^{(\alpha_n-
1)/2}}{n^{\alpha_n}}
\end{equation}
where $\Gamma$ is the Gamma function.
\end{lem}
\proof
Let $F_n(t)=\textrm{Prob}\{\max_i|g_i|\leq t\}$, where $g_1,\dots,g_n$
are independent $N(0,1)$-Gaussian variables. We have
$$\gamma_n(t\cube{n})=\bigg(\sqrt{\haflip}\int_0^te^{-
x^2/2}\,dx\bigg)^n=F_n(t).$$
For $n>1$, put
$$a_n=\frac{1}{\sqrt{2\log n}},\qquad\textrm{and}\qquad b_n=\sqrt{2\log
n}-\frac{\log(\pi\log n)}{2\sqrt{2\log n}}.$$
The well known tail approximation
\begin{equation}\label{eq:tail}\sqrt{\haflip}\int_t^{\8}e^{-x^2/2}\,dx=\sqrt{\haflip}\,
\frac{1+\littleoh{1}}{t}e^{-t^2/2}\qquad\textrm{as}\ \ t\to\8, \end{equation}
combined with a simple calculation, implies that \begin{equation}\label{eq:limit}
\lim_{n\to\8}F_n(a_nx+b_n)=\exp(-e^{-x}),\qquad \lechol x\in\aR. \end{equation}
A change of variables gives:
$$
\aligned
&\int\limits_0^{\8}e^{-\alpha_n
t^2/2}\gamma_n(t\cube{n})\,dt=a_n\int\limits_{-b_n/a_n}^{\8}e^{-
\alpha_n(a_nx+b_n)^2/2}F_n(a_nx+b_n)\,dx\\
&=\frac{\pi^{\alpha_n/2}}{\sqrt{2}}\frac{(\log
n)^{(\alpha_n-1)/2}}{n^{\alpha_n}}e^{- \littleoh{1}}
\int\limits_{-\8}^{\8}e^{-x^2\,\littleoh{1}}e^{-\alpha_n x(1-
\littleoh{1})}F_n(a_nx+b_n)\chi_n(x)
\,dx.\endaligned$$
 Here $\chi_n$ stands for the
characteristic function of the interval
$[-b_n/a_n,\8)$. All four terms of the integrand in the last integral are
non-negative for each $x$. For $x\geq 0$ and sufficiently large~$n$ we have
${e^{-\alpha_n x(1-o(1))}< e^{-\alpha x/2}}$, while the rest of the terms are
majorized by~$1$. For $x<0$ and sufficiently large
$n$, we have ${F_n(a_nx+b_n)<2\exp(-e^{|x|})}$ and ${e^{-\alpha_n
x(1-o(1))}<e^{2\alpha|x|}}$. Thus in both cases if $n$ is sufficiently large, the
integrand is dominated by an integrable function. By (\ref{eq:limit}), the
integrand converges pointwise to the function $e^{-\alpha x}\exp(-e^{-x})$;
Lebesgue's bounded convergence theorem can be applied: \goodbreak
$$\aligned
\lim_{n\to\8}\int_{-b_n/a_n}e^{-x^2\littleoh{1}}e^{-\alpha_n
x(1-\littleoh{1})}F_n(a_nx+b_n)\,dx&=\int_{-\8}^{\8}e^{-\alpha
x}\exp(-e^{-x})\,dx\\
&=\Gamma(\alpha).
\endaligned$$
The proof of Lemma~$2.3$ is
complete.\sofproof
\bigskip
Taking $\alpha_n\equiv k$ in Lemma {\thelem} and bearing in mind
Proposition~$2.2$
re-proves the  following result, which was mentioned in the
introduction.
\begin{cor}
For fixed $k$,
$$f(0,k,n)\sim \frac{2^k}{\sqrt{k}}(\pi\log n)^{(k-1)/2},\qquad{\rm as}\ \ 
n\to\8.$$
\end{cor}
\bigskip
We turn now to the case of fixed co-dimension. The next result is deduced from 
proposition~$2.2$.
\begin{prop}
For $d\geq 1$,
$$f(0,n-d,n)\geq{n\choose d}2^n\left(\frac{1}{\pi}\arctan\frac{1}{\sqrt{n-
d}}\right)^d,\qquad (d\geq 1).$$
Equality holds for $d=1\,:$
\begin{equation}\label{eq:facets2}
f(0,n-1,n)=\frac{2^nn}{\pi}\arctan\frac{1}{\sqrt{n-1}}.
\end{equation}
\end{prop}
\proof  Consider the probability measure $d\mu(t)=2\sqrt{\frac{k}{\pi}}e^{-kt^2}dt$ 
on the
half-line $[0,\8)$. Put
$$\Phi(t)=\frac{2}{\sqrt{\pi}}\int_0^te^{-x^2}\,dx.$$
Then
$$\gamma_{n-k}(t\cube{n-k})=\left(\sqrt{\haflip}\int_0^te^{-
x^2/2}\,dx\right)^{n-k}=\Phi^{n-k}(t/\sqrt{2}).$$
Therefore
\begin{equation}
\aligned
\sqrt{\frac{2k}{\pi}}\int_0^{\8}e^{-kt^2/2}\gamma_{n-k}(t\cube{n-k})\,dt
&=\sqrt{\frac{2k}{\pi}}\int_0^{\8}e^{-kt^2/2}\Phi^{n-k}(t/\sqrt{2})\,dt\\
&=\int_0^{\8}\Phi^{n-k}(t)\,d\mu(t)\\
&\geq\left(\int_0^{\8}\Phi(t)\,d\mu(t)\right)^{n-k}.\endaligned
\end{equation}
Elementary calculation shows that
\begin{displaymath}
\int_0^{\8}e^{-kt^2}\Phi(t)\,dt=\frac{1}{\sqrt{\pi k}}\arctan\frac{1}{\sqrt{k}}.
\end{displaymath}
A combination of ({\theequation}) with proposition $2.2$ gives the desired inequality,
after a replacement of $k$ by $n-d$. Observe that for
$k=n-1$ (that is, $d=1$), there is equality in the inequality of 
({\theequation}).\sofproof\bigskip
\noindent{\bf Remarks}\bigskip
1. For $n=3$ we get from (\ref{eq:facets2}): 
$f(0,2,3)=(24/\pi)\arctan\frac{1}{\sqrt{2}}\approx 4.7$. 
Therefore a random $2$-section of the $3$-cube is more likely to be a parallelogram 
than a
hexagon. 

2. B\'ar\'any and Lov\'asz proved in \cite{BarLov} that (in particular) almost every
$k$-section of the
$n$-cube has at least $2^k$ vertices. Clearly this is a precise lower bound. For $k=n-1$, our
result shows that the expected value is asymptotically $\sqrt{n}/\pi$ times the minimal value.

3. The asymptotic behavior of the integral
$$\int_0^{\8}e^{-kt^2/2}\gamma_{n-k}(t\cube{n-k})\,dt$$
for fixed $k$ and $n\to\8$ was determined in \cite{BorHenk} (following  \cite{Raynaud}),
and was used to prove formula (\ref{eq:BHasymp1}) of the introduction.
See also \cite{AffSch}. The asymptotic result is basically a corollary of the
classical tail approximation of
a single $N(0,1)$-Gaussian variable. Our approach to
the proof of Lemma~$2.3$ seems to simplify the analysis.

4. As was indicated in the introduction, we can choose
$\eps=\frac{c}{\sqrt{\log n}}$ for some constant $c>0$, and then with high
probability a random
$2$-section of the cube is
$(1+\frac{c}{\sqrt{\log n}})$-Euclidean. It is well known that among all centrally symmetric polygons having $2m$ 
vertices, the regular $2m$-gon minimizes the Banach-Mazur distance to the Euclidean disc; the
minimal distance is $(\cos(\pi/2m))^{-1}$. Consequently with high 
probability we have
$$(\cos(\pi/2m))^{-1}<1+\frac{c}{\sqrt{\log n}}.$$ Hence 
 most $2$-sections of the $n$-cube have
at least $C(\log n)^{1/4}$ vertices, for  some
positive constant $C$.  By Corollary $2.2$ (after a suitable rearrangement)
$$f(0,2,n)=2\sqrt{\pi}\ex(\max_{1\leq i\leq n}|g_i|),$$
which is of the order of magnitude of $\sqrt{\log n}$. Summarizing these observations, we
conclude: a typical $2$-section of the $n$-cube is $(1+\frac{c}{\sqrt{\log n}})$-
Euclidean,
hence it cannot have too few vertices --- it has at least $C(\log n)^{1/4}$ vertices 
with
probability that tends to~$1$ as $n\to\8$. It does not however tend to be a regular
polygon, because the expected number of its vertices is too high for that. 

\section{Other faces}
We now turn to the case $j\geq 1$, and prove the following result.
\begin{prop}
For $j\geq 1$, the following inequality holds.
$$
\frac{f(0,k-j,n)}{f(j,k,n)}<
\sqrt{\frac{2}{\pi}}\Big(\frac{j(k-j)}{n-k+j}\Big)^{1/2}
\int_0^{\8}e^{-\frac{j(k-j)}{n-k+j}t^2/2}\gamma_j(t\cube{j})\,dt.
$$
\end{prop}
The starting point in the proof of Proposition $3.1$ is (\ref{eq:basic}) of
Lemma~$2.1$. Again, we 
choose
${A=C_1(\face)}$. The random variable ${g\to\sigma_j(g\sph{k-1}\cap A)}$, which is 
defined on
$\bigoh{n}$, has values in $[0,1]$. Hence
\begin{equation}
\int_{{\bigoh{n}}}\sigma_j(g\sph{k-1}\cap A)\,dH(g)=\int_0^1 H\{g:
\sigma_j(g\sph{k-1}\cap A)\geq t\}\,dt.
\end{equation} 
The integrand is non-increasing, and 
\begin{equation}
H\{g: \sigma_j(g\sph{k-1}\cap A)\geq
0\}=\textrm{Prob\,}\{X\cap\face\neq\emptyset\},
\end{equation}
 because the event
$\{g\sph{k-1}\cap A\neq\emptyset\enspace\textrm{and}\enspace\sigma_j(g\sph{k-
1}\cap
A)=0\}$ has Haar measure zero. Therefore by (\ref{eq:basic}):
$$\aligned
\sigma_{n-k+j}(A)&\leq \textrm{Prob\,}\{X\cap\face\neq\emptyset\}\sup\{t: H\{g:
\sigma_j(g\sph{k-1}\cap A)\geq t\}>0\}\cr &\leq
\textrm{Prob\,}\{X\cap\face\neq\emptyset\}\sup\{\sigma_j(g\sph{k-1}\cap A): g\in
\bigoh{n}\}.
\endaligned$$
Let 
$$
t_{j,k,n}=\sup\{\sigma_j(g\sph{k-1}\cap A): g\in\bigoh{n}\}.
$$
By (\ref{eq:basic}), $\numbeforethat$ and $\numbefore$ we get
$$\textrm{Prob}\{X\cap\face\neq\emptyset\}\geq\frac{\sigma_{n-
k+j}(A)}{t_{j,k,n}}.$$
Hence by (\ref{eq:Prob})
$$f(j,k,n)\geq 2^{k-j}{n\choose
k-j}\frac{\sigma_{n-k+j}(A)}{t_{j,k,n}}=\frac{\frac{1}{2}f(0,k-
j,n)}{t_{j,k,n}}.$$
We must bound $t_{j,k,n}$ from above.   Since
$A$ is contained in a half-space, a trivial bound is $t_{j,k,n}\leq\frac{1}{2}$. In 
some
cases this bound can be significantly improved. The main lemma in this section
is the following.
\begin{lem}\label{lem:bound}
If $1\leq j<k<n$, then
$$t_{j,k,n}\leq
\frac{1}{\sqrt{2\pi}}\Big(\frac{j(k-j)}{n-k+j}\Big)^{1/2}
\int_0^{\8}e^{-\frac{j(k-j)}{n-k+j}t^2/2}\gamma_j(t\cube{j})\,dt.
$$
\end{lem}

The next lemma will be used in the proof of Lemma~\ref{lem:bound}.
\begin{lem}
Given a positive number $\tau>0$, a $j$-dimensional subspace $Y$ of $\aR^m$ and
a point $y_0\in Y$, the following inequality holds.
\begin{equation}\label{eq:braslieb}
\gamma_j((Y\cap\tau\cube{m})-y_0)\leq\gamma_j(\tau\sqrt{m/j}\cube{j}).
\end{equation}
\end{lem}
\vskip 6pt
\proof 
Let $Q$ denote the orthogonal projection onto~$Y-y_0$. As usual,
$\{e_i\}_{i=1}^m$ are
the standard unit vectors in $\aR^m$. Put ${u_i=Qe_i/\norma{Q e_i}}$ if $Q e_i\neq 0$, and $u_i=0$ otherwise; put 
$c_i=\norma{Q
e_i}^2$ and $\alpha_i=\ip{y_0}{e_i}$ for $1\leq i\leq
m$. ($\ip{\cdot}{\cdot}$ is the standard scalar product.) Then
$$\aligned
Y\cap\tau\cube{m}&=\{y\in Y: |\ip{y}{e_i}|\leq\tau\ \ \lechol i\}\cr
&=\{y\in Y: |\ip{y-y_0}{e_i}+\ip{y_0}{e_i}|\leq\tau\ \ \lechol i\}\cr
&=\{y\in Y:
\frac{-\alpha_i-\tau}{\sqrt{c_i}}\leq\ip{y-
y_0}{u_i}\leq\frac{-\alpha_i+\tau}{\sqrt{c_i}}\}.
\endaligned$$
Therefore
$$(Y\cap\tau\cube{m})-y_0=\{x\in Y-y_0:
\frac{-\alpha_i-\tau}{\sqrt{c_i}}\leq\ip{x}{u_i}\leq\frac{-\alpha_i+\tau}{\sqrt{c_i}}\}$$
Now we can imitate K.~Ball's argument from \cite{Balcub} concerning sections of maximal
volume. Instead of the Lebesgue measure, we have to consider the Gaussian measure.  

In $Y-y_0$, the identity operator can be written as $\sum_1^{m}c_iu_i\otimes 
u_i$.
In particular,
$$\sum_{i=1}^{m}c_i=j,\quad\textrm{and}\quad\norma{x}^2=\sum_{i=1}^{m}c_i\ip{x}{u_i}^2,\quad\forall
x\in Y-y_0.$$  
Therefore the Gaussian measure in $Y-y_0$ is equal to
$$(2\pi)^{-j/2}\exp(-\sum_{i=1}^{m}c_i\ip{x}{u_i}^2/2)\,dx.$$
Let $\chi_i$ denote the characteristic function of the
interval $[\frac{-\alpha_i-\tau}{\sqrt{c_i}},\frac{-\alpha_i+\tau}{\sqrt{c_i}}]$. Then, by
the above,
\begin{equation}\aligned
\gamma_j((Y\cap\tau\cube{m})-y_0)&=(2\pi)^{-j/2}\int_{Y-
y_0}\Big(\prod_{i=1}^{m}
\chi_i(\ip{x}{u_i})e^{-c_i\ip{x}{u_i}^2/2}\Big)\,dx\cr
&=(2\pi)^{-j/2}\int_{Y-y_0}\prod_{i=1}^{m}(\chi_i(\ip{x}{u_i})e^{-
\ip{x}{u_i}^2/2})^{c_i}\,dx\cr
&\leq(2\pi)^{-j/2}\prod_{i=1}^{m}\big(\int_{(-\alpha_i-\tau)/\sqrt{c_i}}^{(-\alpha_i+\tau)/\sqrt{c_i}}e^{-s^2/2}\,ds\big)^{c_i}.\cr
\endaligned\end{equation}
The last inequality is a consequence of Brascamp-Lieb's inequality, which is stated in
\cite{Balcub} as follows: \par
{\bf Lemma}\enspace {\it Let $(u_i)_1^m$ be a sequence of unit vectors in $\aR^n$ and $(c_i)_1^m$ a sequence
of positive numbers so that
$$\sum_1^mc_iu_i\otimes u_i=I_n.$$
For each $i$, let $f_i:\aR\to [0,\8)$ be integrable. Then
$$\int_{\aR^n}\prod_{i=1}^mf_i(\ip{u_i}{x})^{c_i}\,dx\leq\prod_{i=1}^m\left(\int_{\aR}f_i\right)^{c_i}.$$
}
The $i$'th integral in the product of $\numnow$ is not larger than
$\int_{-\tau/\sqrt{c_i}}^{\tau/\sqrt{c_i}}e^{-s^2/2}\,ds$.
Hence the last expression in $\numnow$ is bounded above by
$$(2\pi)^{-j/2}\prod_{i=1}^{m}\big(2\int_0^{\tau/\sqrt{c_i}}e^{-
s^2/2}\,ds\big)^{c_i},$$
which is maximized when all the $c_i$'s are equal. Hence
\goodbreak
\begin{equation}
\aligned
\gamma_j((Y\cap\tau\cube{m})-y_0)&\leq
\big(\sqrt{\frac{2}{\pi}}\int_0^{\tau\sqrt{m/j}}e^{-s^2/2}\,ds\big)^j\cr
&=\gamma_j(\tau\sqrt{m/j}\cube{j}).\endaligned
\end{equation}
The proof of Lemma $3.3$ is complete.\sofproof
\bigskip
\noindent\textbf{Proof of Lemma \ref{lem:bound}}
\vskip 6pt
For $g\in\bigoh{n}$
$$\aligned
\sigma_j(g\sph{k-1}\cap
A)&=\gamma_{j+1}\big(C(\face)\cap\textrm{span}(g\sph{k-1})\big)\cr
&=\gamma_{j+1}\big(C[\face\cap\textrm{span}(g\sph{k-1})]\big).\endaligned$$
The second equality is a consequence of the identity $C(\face)\cap X= C(\face\cap X)$,
which trivially holds for every subspace $X\subset\aR^n$. 
Fix a subspace $X\in G_{n,k}$ for which the section ${X\cap\face}$ is 
$j$-dimensional;
almost every $X\in G_{n,k}$ has this property. Let $C$ denote the
$(j+1)$-dimensional cone generated by $X\cap\face$; put $X_0=\textrm{span}C$.  By 
$M$
we denote the affine subspace spanned by~$X\cap\face$, and by
$d$, its distance from the origin of~$X$. The Gaussian measure of the cone
$C$ is computed as follows. Take the unit vector
$\xi\in X_0$ which is orthogonal to~$M$, and for which $d\xi\in M$. For $t>0$, put
${W_t=\{x\in X_0:\ip{x}{\xi}=t\}}$.
Observe that
$C\cap W_t=(t/d)(X\cap\face)$.
Let $P$ denote the orthogonal projection from $X_0$ onto $W_0$.  By Fubini's theorem:
\begin{equation}\label{eq:conemes}
\aligned
\gamma_{j+1}(C)&=\frac{1}{\sqrt{2\pi}}\int_0^{\8}e^{-t^2/2}
\gamma_j(P(C\cap W_t))\,dt\cr
&=\frac{1}{\sqrt{2\pi}}\int_0^{\8}e^{-t^2/2}\gamma_j(P(t/d)(X\cap\face))\,dt.
\endaligned
\end{equation}
Our task is to estimate the expression $\gamma_j(P\tau(X\cap\face))$
for every $\tau>0$.
We will need to discuss Gaussian measures in different subspaces. Whenever $M$ is an $m$-dimensional 
subspace of~$\aR^n$ and ${q\in M}$, let $\gee_{M,q}$ denote the measure
$(2\pi)^{-m/2}\exp(-\norma{x-q}^2/2)\,dx$. In
case $M$ is an $m$-dimensional linear subspace of $\aR^n$ and $q=0$ we shall simply write 
$\gee_{M,0}=\gamma_m$. If $T$
is an isometry of $\aR^n$, then for every Borel subset $S\subset M$ we have
\begin{equation}\label{eq:gsmeas}
\gee_{M,q}(S)=\gee_{TM,Tq}(TS).
\end{equation}
Let us momentarily assume  that
$\tau=1$.  Let
$q$ denote the nearest point of~$M$ to the origin of~$X$. Both $M$ and the range of the projection~$P$
are $j$-dimensional affine subspaces of~$X_0$. We have
$$P(X\cap\face)=(X\cap\face)-q,$$
hence by (\ref{eq:gsmeas})
\begin{equation}
\gee_{M,q}(X\cap\face)=\gamma_j(P(X\cap\face)).
\end{equation}
Now let $L$ denote the affine subspace
spanned by
$\face$, whose origin $O_L$ is taken as the center of the face $\face$. (So if $X$ passes
through the center of $\face$, then $q=O_L$.) $M$ is also a $j$-dimensional
affine subspace of~$L$.
By~(\ref{eq:gsmeas}),
$$\gee_{M,q}(X\cap\face)=\gee_{M-(q-O_L),O_L}\big((X\cap\face)-(q-O_L)\big).$$
Applying the same argument for arbitrary
$\tau>0$ we conclude that
\begin{equation}\label{eq:projeq}
\gamma_j(P\tau(X\cap\face))=\gee_{\tau M-\tau(q-O_L),\tau O_L}\big(\tau(X\cap\face)-\tau(q-O_L)\big).
\end{equation}

We may think of $\tau L$ as $\aR^{n-k+j}$, of $\tau\face$ as $\tau\cube{n-k+j}$, and of
$\tau(X\cap\face)$ as an affine $j$-dimensional section of $\tau\cube{n-k+j}$.
Thus for each $t>0$ Lemma $3.3$ can be used with $\tau=t/d$ and ${m=n-k+j}$. 
By the definition of~$d$,
we have
$d\geq\sqrt{k-j}$. Combining (\ref{eq:braslieb}),(\ref{eq:conemes}) and $\numnow$ we deduce that
\goodbreak
$$\aligned
\gamma_{j+1}(C)&\leq\frac{1}{\sqrt{2\pi}}\int_0^{\8}e^{-
t^2/2}\gamma_j\bigg(t\Big(\frac{n-k+j}{j(k-j)}\Big)^{1/2}\cube{j}\bigg)\,dt\cr
&=\frac{1}{\sqrt{2\pi}}\Big(\frac{j(k-j)}{n-k+j}\Big)^{1/2}\int_0^{\8}
\exp(-\frac{j(k-j)}{n-k+j}t^2/2)\gamma_j(t\cube{j})\,dt.\cr\endaligned
$$
The proof of lemma~\ref{lem:bound} and thus of proposition~$3.1$ is complete.
\sofproof
\bigskip
By using the asymptotic formulas of section~$2$, namely Lemma~$2.3$ and
Corollary~$2.4$, we can now prove the following result, which shows that the lower bound
for $f(j,k,n)$ derived from proposition~$3.1$ is, in some cases, asymptotically best possible.
\begin{cor}
For fixed integers $1\leq l<m$,
\begin{equation}
f(n-m,n-l,n)\backsim \frac{(2n)^{m-l}}{(m-l)!}\quad\textrm{as}\ \ n\to\8.
\end{equation}
\end{cor}
\proof  Put $\alpha_n=(m-l)(n-m)/(n-m+l)$. By Proposition $3.1$,
\begin{equation}
\frac{f(0,m-l,n)}{f(n-m,n-l,n)}<\sqrt{\frac{2\alpha_n}{\pi}}\int_0^{\8}e^{-\alpha_nt^2/2}
\gamma_{n-m}(t\cube{n-m})\,dt.
\end{equation}
Put $b_n=(\log(n-m))^{(\alpha_n-1)/2}/(n-m)^{\alpha_n}$ and $c_n=(\log n)^{(m-l-1)/2}$. 
Let $d_n$ denote the right hand side of~ $\numnow$, from which
we get
$$f(n-m,n-l,n)\frac{b_n}{c_n}>\frac{f(0,m-l,n)}{c_n}\frac{b_n}{d_n}.$$
Since $\lim_{n\to\8}\alpha_n=(m-l)$, Lemma~$2.3$ implies that
\begin{displaymath}
\lim_{n\to\8}\frac{b_n}{d_n}=\frac{1}{\pi^{(m-l-1)/2}\Gamma(m-l)\sqrt{m-l}}.
\end{displaymath}

Moreover, by Corollary $2.4$,
\begin{displaymath}
\lim_{n\to\8}\frac{f(0,m-l,n)}{c_n}=\frac{2^{m-l}\pi^{(m-l-1)/2}}{\sqrt{m-l}}.
\end{displaymath}
Thus, the sequence $f(n-m,n-l,n)\frac{b_n}{c_n}$ is larger than a sequence that
tends to $2^{m-l}/(m-l)!$ as $n$ tends to infinity. On the other hand
we have ${f(n-m,n-l,n)<2^{m-l}{n\choose m-l}}$, so
\goodbreak
$$f(n-m,n-l,n)\frac{b_n}{c_n}<2^{m-l}{n\choose m-l}\frac{b_n}{c_n},$$
and since $b_n/c_n\backsim n^{l-m}$, the r.h.s here tends to $2^{m-l}/(m-l)!$. Consequently,
$$\lim_{n\to\8}f(n-m,n-l,n)\frac{b_n}{c_n}=\frac{2^{m-l}}{(m-l)!}.$$
The required asymptotic formula  
follows immediately. The proof of Corollary~$3.4$ is complete.\sofproof
\bigskip\goodbreak
\noindent{\bf Remarks}
\bigskip
1. The previous corollary implies that the number of $(n-m)$-faces of a random $(n-l)$-section of the $n$-cube tends to concentrate near the value $2^{m-l}{n\choose m-l}$, which bounds it from above. So
for example, a typical $1$-co-dimensional section of the $n$-cube will have $2n-o(n)$ facets as
$n\to\8$. This result can also be deduced from the identity (\ref{eq:AffSch}). Indeed, by
duality we have $f(n-m,n-l,n)=\ex(f_{m-l-1}(\Pi_{n-l}(B_1^n)))$, and replacing $T^n$ by $B_1^n$
in the proof of Theorem~2 in \cite{AffSch}, (the details of this replacement appear in \cite{BorHenk};
see the proof of Theorem~1.1 there) we get the previous corollary.

2. 
According to a remark made in \cite{BorHenk}, the 
number
$f(j,k,n)$ is equal to the expected number of $(k-j-1)$-faces of the convex hull of 
$\pm
G_1,\dots,\pm G_n$,  where the $G_i$'s are independent copies of a $k$-dimensional 
Gaussian
vector. Hence, 
the
results for $f(0,k,n)$ can be interpreted as results for the expected number of facets of
the convex hull of $\{\pm G_i \}_1^n$ in $\aR^{k}$. For example, we can 
translate the first remark at the end of section~$2$ to the following statement:
\vskip 6pt
{\it If\, $3$ points in the plane are chosen at random, then their symmetric convex hull is more likely to be a parallelogram than a hexagon\/}.
\vskip 6pt
\bigskip
\noindent\textbf{Acknowledgements}\quad I thank Itai Benjamini for getting me interested 
in
the subject of this paper, and for several stimulating discussions. I also thank the referees
for their useful comments.
\bibliographystyle{amsplain}
\providecommand{\bysame}{\leavevmode\hbox to3em{\hrulefill}\thinspace}

\bigskip
\begin{tabular}{ll}
Yossi Lonke & Current address:\\
Institute of Mathematics & Department of Mathematics \\
The Hebrew University of Jerusalem & Case Western Reserve University \\
Jerusalem 91904, Israel & Cleveland, Ohio 44106-7058\\
{} &email: {\emailfont jrl16@po.cwru.edu}\\
\end{tabular}
\end{document}